\documentclass{article}%
\usepackage{amsfonts}
\usepackage{amsmath}
\usepackage{amssymb}
\usepackage{graphicx}%
\def\thebibliography#1{\small                                        \list{\arabic{enumi}.}                                         {\settowidth\labelwidth{[#1]}                                 \leftmargin\labelwidth                                        \itemsep 0pt                                                  \parsep \itemsep                                              \advance\leftmargin\labelsep                                  \usecounter{enumi}}}
\def\@cite#1{$^{#1}$}

\newenvironment{myabstract}{\normalsize\noindent}{}

\begin{document}

%*
%*
%%% The title of your contribution
%%% please do not force line breaks with \\
%%% even if you think it looks ugly

\setcounter{equation}{0}\setcounter{footnote}{0}\begin{flushleft}{\bf%
\textbf{\Large Studies on concave Young-functions}}\end{flushleft}\vspace
{0.25cm}

%%%
%%% Your name and affiliation
%%%

\leftline{{\large         N. Kwami Agbeko}}\vspace{-0.05cm}

{\small \begin{flushleft}{\it  Institute of Mathematics, University of
Miskolc, 3515 Miskolc-Egyetemvaros, Hungary}\\{\tt
email: matagbek@uni-miskolc.hu}\\{\tt
Subject Classification: Primary 47H10, 47H11, 26A18; Secondary 26A06, 26A09, 33B30, 37C25}
\end{flushleft}
%This line gives your e-mail address
%If you don't want it to appear, please
%type {} instead
}

{\small
%%% Each paper should contain the 2000 Mathematics Subject Classification.
}

\bigskip

\begin{center}
\textbf{ABSTRACT}
\end{center}

\begin{myabstract}
We succeeded to isolate a special class of concave Young-functions enjoying
the so-called \emph{density-level property}. In this class there is a proper
subset whose members have each the so-called degree of contraction denoted by
$c^{\ast}$, and map bijectively the interval $\left[  c^{\ast},\,\infty
\right)  $ onto itself. We constructed the fixed point of each of these
functions. Later we proved that every positive number $b$ is the fixed point
of a concave Young-function having $b$ as degree of contraction. We showed
that every concave Young-function is square integrable with respect to a
specific Lebesgue measure. We also proved that the concave Young-functions
possessing the density-level property constitute a dense set in the space of
concave Young-functions with respect to the distance induced by the $L^{2}$-norm.
\end{myabstract}

We note that the full paper published in {\small \emph{Miskolc Mathematical
Notes}} can be downloaded at the homepage of the above periodical here
inclosed and that Lemma 7 is wrong. Nevertheless, nothing is lost cf. the
referee's note for the \emph{Mathematical Review}:
\textbf{MR2148839(2006e:26005)}

\end{document}